\documentclass[brochure,english,12pt]{smfbourbaki}

\usepackage[T1]{fontenc}
\usepackage{lmodern,amssymb,bm}
\usepackage{mathrsfs}
\usepackage{babel}

\usepackage[utf8]{inputenc}

\usepackage[colorlinks=true, linkcolor=blue, citecolor=red, urlcolor=blue]{hyperref}

\usepackage[
backend=bibtex,
style=authoryear, 
citestyle=authoryear-comp,
maxnames=7,
sortcites=false 
]{biblatex}
\usepackage{csquotes}

\newtoggle{tnbcbx@howcited}
\DeclareEntryOption[boolean]{howcited}[true]{\settoggle{tnbcbx@howcited}{#1}}
\DeclareBibliographyOption[boolean]{howcited}[true]{\settoggle{tnbcbx@howcited}{#1}}
\DeclareTypeOption[boolean]{howcited}[true]{\settoggle{tnbcbx@howcited}{#1}}

\newbibmacro{howcited}{%
  \iftoggle{tnbcbx@howcited}
    {\iffieldundef{shorthand}
       {}
       {\setunit{\addspace}%
        \printtext[parens]{%
          \bibstring{citedas}%
          \setunit{\addcolon\space}%
          \printfield{shorthand}%
          \setunit{\addslash}%
          }}}
    {}}

\renewbibmacro{finentry}{\usebibmacro{howcited}\finentry}

\DefineBibliographyStrings{french}{citedas    = {cité ci-dessus avec l'acronyme}}
\DefineBibliographyStrings{english}{citedas    = {heretofore cited as}}

\DefineBibliographyExtras{french}{\restorecommand\mkbibnamefamily}

\DeclareDelimFormat{nameyeardelim}{\addcomma\space}

\DeclareNameAlias{sortname}{given-family}

\renewcommand{\bibnamedash}{\leavevmode\raise3pt\hbox to3em{\hrulefill}\space}

\AtEveryBibitem{%
  \clearfield{issn} 
  \clearfield{isbn} 
  \clearfield{doi} 
  \clearlist{language} 
  \ifentrytype{online}{}{
  \ifentrytype{unpublished}{}{
    \clearfield{url}
  }
  }
}

\renewbibmacro{in:}{%
    \ifentrytype{article}{}{\printtext{\bibstring{in}\intitlepunct}}}

\DeclareFieldFormat[article,periodical,inreference]{number}{\mkbibparens{#1}}
\DeclareFieldFormat[article,periodical,inreference]{volume}{\mkbibbold{#1}}
\renewbibmacro*{volume+number+eid}{%
    \printfield{volume}%
    \setunit*{\addthinspace}
    \printfield{number}%
    \setunit{\addcomma\space}%
    \printfield{eid}}

\DeclareFieldFormat[article,inbook,incollection]{title}{\enquote{#1}\addcomma} 

\date{Janvier 2024}
\bbkannee{76\textsuperscript{e} année, 2023--2024}  
\bbknumero{1216}                                      

\title{Exotic subgroups of hyperbolic groups}

\author{Olivier Guichard}
\address{IRMA, Université de Strasbourg\\
67000 Strasbourg}
\email{guichard@unistra.fr}

\newcommand{\FH}{\mathrm{FH}}
\newcommand{\FP}{\mathrm{FP}}

\begin{document}

\maketitle

\section{Introduction}
\label{sec:introduction}
Since their introduction by \textcite{gromov}, word hyperbolic groups have been the
focus of a lot of activity and have proved central in attacking a number of
problems. It was soon noticed that their cohomological
properties are very strong. As a matter of fact, a torsion-free word hyperbolic
group~$\Gamma$ is of type~$\mathcal{F}$, meaning that $\Gamma$~is the
fundamental group of a finite aspherical cell complex (\textcite{gromov} attributes this to Eliyahu Rips). Such a property for a
group~$\Gamma$ is called a \emph{finiteness} property. For every positive
integer~$n$, there is a coarser finiteness property
denoted~$\mathcal{F}_n$ that requires a group~$\Gamma$ to be the fundamental
group of an aspherical cell complex, possibly infinite, but which has only finitely many
cells up to dimension~$n$. A group of type~$\mathcal{F}_n$ and not of
type~$\mathcal{F}_{n+1}$ is sometimes said to have exotic finiteness
properties\footnote{This terminology was coined down by
  \textcite[section~5]{les-roumains}, later used in a book review by
  \textcite{meier}, and popularized by \textcite{fourier}; see also 
 the
title of section~7 in \textcite{jankiewicz-norin-wise}. A formal definition
appeared first in \citeauthor{isenrich-py-exotic} (\citeyear*{isenrich-py-exotic}).}. The aim
of this report is to illustrate that word hyperbolic groups can have exotic
subgroups: subgroups with exotic finiteness properties or subgroups of
type~$\mathcal{F}$ but not word hyperbolic.

\begin{theo}[{\citeauthor{llosa-pi-inventiones}, \citeyear*{llosa-pi-inventiones}, corollary~3}]
  \label{theo:llosa_py}
  Let~$n$ be a positive integer. There exists a word hyperbolic group~$\Gamma$
  containing a subgroup that is of type~$\mathcal{F}_n$ but not of type~$\mathcal{F}_{n+1}$.
\end{theo}

\begin{theo}[{\citeauthor{italiano-martelli-migliorini5},
 \citeyear*{italiano-martelli-migliorini5}, corollary~2}]
  \label{theo:les_italiens}
  There exists a word hyperbolic group~$\Gamma$ containing a subgroup of
  type~$\mathcal{F}$ that is not word hyperbolic.
\end{theo}

In both statements the subgroups are  kernels of  homomorphisms
from~$\Gamma$ to~$\mathbf{Z}$ (and in particular are normal subgroups). The
geometric counterparts of these homomorphisms are maps from~$M$ to the circle,
where $M$~is a manifold or a pseudo-manifold whose fundamental group is~$\Gamma$. This makes the
analysis of the subgroups amenable to Morse-theoretical techniques. More
precisely, on one hand Lefschetz theory is used by
\citeauthor{llosa-pi-inventiones} (\citeyear*{llosa-pi-inventiones})
to study word hyperbolic groups 
which are arithmetic subgroups of
$\mathrm{U}(n,1)$ (and $M$~is the quotient of the unit ball in~$\mathbf{C}^n$
by the action of these arithmetic subgroups); on the other hand the Morse theory of affine
cell complexes developed by \textcite{bestvina-brady} is used by
\citeauthor{italiano-martelli-migliorini5}
 (\citeyear*{italiano-martelli-migliorini5}) for their word hyperbolic groups which are given by combinatorio-geometrical data.

As emphasized by the authors themselves (and apparent in that the
above citations point to corollaries), the interesting statements may not be
the above results, that give positive solutions to questions raised
after the introduction of word hyperbolic groups, but the geometric
constructions of which they are the shadows. The present report will indeed
sketch these constructions and try to refer to the original articles for complete proofs.

\subsection*{Acknowledgment}
\label{sec:acknowledgment}

I am grateful to Nicolas Bourbaki for offering the opportunity to present this
seminar and for their careful readings of this text. I am glad to thank here Claudio Llosa Isenrich, Bruno Martelli, and Pierre Py for
their help while writing these notes. The feedback of Claudio and Pierre
has improved this text far beyond what the author could have written on his own.

\section{A brief overview of the historical development and further statements}
\label{sec:brief-overv-hist}

Finiteness properties have various declinations: 
the properties $\FH_n(R)$ ($R$~is a ring) request that
the group~$\Gamma$ is the fundamental group of a compact cell complex whose universal cover
has trivial reduced homology
with coefficients in~$R$ in degrees $<n$ \parencite[pp.~445--446]{bestvina-brady}, and the properties $\FP_n(R)$ request that
the trivial $R[\Gamma]$-module has a projective resolution whose
homogeneous factors of degree $\leq n$ are
finitely generated. The property~$\mathcal{F}_1$ is equivalent to the
group~$\Gamma$ being finitely generated.  The property~$\mathcal{F}_2$ is equivalent to the
group~$\Gamma$ being finitely presented. For all~$n$, the
property~$\mathcal{F}_n$ implies the property $\FH_n(\mathbf{Z})$,
$\FH_n(\mathbf{Z})$ implies $\FH_n({R})$, and $\FH_n({R})$
implies $\FP_n({R})$ (and $\FP_n(\mathbf{\mathbf{Z}})$ implies $\FP_n({R})$).

\textcite[corollary (b)]{rips} constructed the first example of a finitely
generated, hence~$\mathcal{F}_1$, but not finitely presented, hence
not~$\mathcal{F}_2$, subgroup in a small cancellation group (in particular in
a word hyperbolic group).
In his essay, \textcite[section 4.4.A]{gromov} suggested a strategy for finding
subgroups with exotic finiteness properties 
  in a word hyperbolic group, by taking
 covers of a flat torus, ramified over a union of codimension~$2$ tori meeting
 orthogonally, and that
fiber over the circle 
  (later Mladen Bestvina showed that Gromov's construction does not lead to a word hyperbolic
group; his argument is reproduced in \textcite[pp.~70--71]{BradyRileyShort}).

The question of the existence, in word hyperbolic groups, of subgroups of type
$\mathcal{F}_n$ and not $\mathcal{F}_{n+1}$ was explicitely raised by
\textcite[p.~130]{gersten} (who uses the notation $\FP_n$ instead of the now
established notation~$\mathcal{F}_n$).
It was also
stated by \textcite[question~7.1]{brady} who constructed finitely presented
subgroups (hence of type~$\mathcal{F}_2$) not of type
$\mathcal{F}_3$\footnote{Brady asks the existence of a finitely presented
  subgroup of type $\FP_n(\mathbf{Z})$ and not $\FP_{n+1}(\mathbf{Z})$. However for
a finitely presented group, the implication $\FP_n(\mathbf{Z})\Rightarrow
\mathcal{F}_n$ holds (a proof can be found in the
  proof of theorem~7.1 in \textcite[chapter VIII]{brown}, it relies on the
  Hurewicz theorem), thus Brady's question is indeed Gersten's question.}.
More examples of finitely presented and not $\mathcal{F}_3$ subgroups,
elaborating on Brady's construction and building on the Bestvina--Brady Morse
theory (see below section~\ref{sec:bestvina-brady-morse}), were subsequently obtained by \textcite{Lohda,Kropholler,kropholler2023finitely}.
\textcite{isenrich-martelli-py} built  the first example of a subgroup of
 type~$\mathcal{F}_3$ and not~$\mathcal{F}_4$ elaborating on a fibration of a
 complete, finite volume, hyperbolic $8$-manifold constructed in
 \citeauthor{italiano2022hyperbolic} (\citeyear{italiano2022hyperbolic}) and
 gave examples of subgroups of type
$\FP_n(\mathbf{Q})$ and not $\FP_{n+1}(\mathbf{Q})$ in cubulable arithmetic lattices of
the Lie group~$\mathrm{O}(2n,1)$.

Kernels of homomorphisms onto~$\mathbf{Z}$ give examples of groups with intermediate
finiteness properties. For example the kernel of the morphism from the free
group~$\mathbb{F}_2$ onto~$\mathbf{Z}$ mapping all the generators to~$1$ is not
finitely generated; the kernel of the morphism $\mathbb{F}_2\times
\mathbb{F}_2 \to \mathbf{Z}$ sending every generator to~$1$ is finitely
generated but not finitely presented. \textcite{stallings}
gave the first example of a group of type~$\mathcal{F}_2$ (thus finitely
presented) that is not of type~$\mathcal{F}_3$; it was later observed
\parencite{gersten} that this example is isomorphic to the kernel of the
morphism from $(\mathbb{F}_2)^3$ to~$\mathbf{Z}$ sending every generator
to~$1$. For every positive integer~$n$, the kernel of the similar homomorphism
from $(\mathbb{F}_2)^n$ to~$\mathbf{Z}$
is of type~$\mathcal{F}_{n-1}$ and not of type~$\mathcal{F}_n$
\parencite{bieri76}. 

\smallskip

On the other hand, the question (answered thus negatively by
theorem~\ref{theo:les_italiens}) whether a subgroup of type~$\mathcal{F}$ in a
word hyperbolic group is itself hyperbolic can be traced back to Bestvina's
problem list\footnote{Written in August 2000 and available at \url{https://www.math.utah.edu/~bestvina/},
retrieved on January 12th 2024.} and is also stated by
\textcite[question~7.2]{brady}. More recently the question appears in
\textcite[section~7]{jankiewicz-norin-wise}. The techniques developed in this
previous reference have been used by \citeauthor{italiano-martelli-migliorini5}
(\citeyear{italiano2022hyperbolic,italiano-martelli-migliorini5}) to
construct fibrations of hyperbolic manifolds over the circle and the
fibration of a pseudo-manifold
explained below in section~\ref{sec:constr-from-right} that leads to
theorem~\ref{theo:les_italiens}. Constructions of hyperbolic manifolds
along the same line were also proposed in
\textcite{kolpakov-slavich,kolpakov-martelli}.

The related question whether a finitely presented subgroup of a word
hyperbolic group of cohomological dimension~$2$ is itself hyperbolic has a
positive answer \parencite{Gersten-codim2}. The similar question in
dimension~$3$ or~$4$ (is it true that an $\mathcal{F}_3$, resp.\ $\mathcal{F}_4$, subgroup of a hyperbolic group of
cohomological dimension~$3$, resp.~$4$, is hyperbolic) is still open. In dimension~$5$,
theorem~\ref{theo:les_italiens} provides a counter-example.

\smallskip

The discussion so far emphasizes morphisms onto~$\mathbf{Z}$. Central
objects which we will not discuss, but enable a finer understanding of
the
finiteness properties of the
kernels of these morphisms,  are the Bieri--Neumann--Strebel invariant \parencite{BNS} and its higher
degree relatives introduced by Renz \parencite{bieri-renz,renz-these,renz-singapour} (the BNSR
invariants). \citeauthor{isenrich-py-exotic} (\citeyear{isenrich-py-exotic})
give other constructions of subgroups of Kähler groups with exotic
finiteness properties. Certain constructions use morphisms to higher-rank
Abelian groups and  are
not amenable to the strategy we describe below, but rely on the BNSR
invariants. Theorem~1.4 in the previous reference constructs subgroups of (not
word hyperbolic) Kähler groups with
intermediate finiteness properties that are
not normal and are themselves Kähler (the construction there involves fiber products
rather than morphisms). \textcite{les-roumains} constructed the first examples
of Kähler groups with intermediate finiteness properties and their techniques
(maps to elliptic curves) were pushed further by others; we refer to
\citeauthor{isenrich-py-exotic} (\citeyear{isenrich-py-exotic}, section~3.1)
for a discussion as well as other references.

\smallskip

The $\ell^2$-homology also gives control on the BNSR invariants
and on finiteness properties of kernels. A consequence of a theorem
of \textcite{Luck} implies that the kernel of a surjective morphism $G\to
\mathbf{Z}$ has not type $\FP_n(\mathbf{Q})$ as soon as the $n$-th
$\ell^2$-Betti number of~$G$ is nonzero.
For the class of residually finite rationally solvable groups (cf.\
\cite{Agol}, for a definition),
\textcite[for the case $n=1$]{kiekak-RFRS} and \textcite[for the general
case]{fisher2022improved} proved that the $\ell^2$-Betti numbers of~$G$ vanish
up to degree~$n$ if and only if there is a surjective morphism $G_1\to
\mathbf{Z}$ with kernel of type $\FP_n(\mathbf{Q})$ where $G_1$~is a finite
index subgroup of~$G$. This was involved in the result of
\textcite{isenrich-martelli-py} mentioned above.

\section{A construction from complex geometry}
\label{sec:constr-from-compl}

Hereafter the article \citeauthor{llosa-pi-inventiones}
 (\citeyear*{llosa-pi-inventiones}) will be
mentioned as \cite{llosa-pi-inventiones} and the article  \citeauthor{isenrich-py-exotic}
(\citeyear*{isenrich-py-exotic}) will be
mentioned as \cite{isenrich-py-exotic}.

 In this section we address theorem~\ref{theo:llosa_py}. The construction
here has three steps. First the kernels of rational cohomology classes of degree~$1$ coming
from complex geometry (precisely 
admitting a Morse representative that is the real part of a complex differential form with isolated
zeros on a Kähler manifold) are shown to produce the wanted example. Second
finite-to-one maps to complex tori provide such cohomology classes. Finally
some arithmetic
quotients of the unit ball in~$\mathbf{C}^n$ immerse into their Albanese
varieties and thus admit finite-to-one maps to a complex torus. This is the strategy
developed in \cite{llosa-pi-inventiones} with a simplification suggested in
\cite{isenrich-py-exotic} (section 8) avoiding
the use of the BNSR invariants. 

\subsection{Forms with isolated zeroes}
\label{sec:forms-with-isolated}

Let~$X$ be a compact connected complex manifold. A closed holomorphic
$1$-form~$\alpha$ on~$X$ leads to a real differential form $a=\Re \alpha$ that represents an element in
the first cohomology group $H^1(X; \mathbf{R})$. When this form is rational,
i.e.\ when the class of~$a$ belongs to $H^1(X;\mathbf{Q})= \mathrm{Hom}(\pi_1(X),
\mathbf{Q})$, it gives rise to a homomorphism from $\pi_1(X)$ onto
a finitely generated subgroup of~$\mathbf{Q}$; hence, up to scaling, it is a
surjective homomorphism from~$\pi_1(X)$
onto~$\mathbf{Z}$. When $X$~is aspherical and $\alpha$~has finitely many zeroes, the kernel of this homomorphism has
the desired exotic finiteness properties. 

\begin{prop}[{\cite[theorem 6.(1)]{llosa-pi-inventiones}}]
  \label{prop:forms-with-isolated}
  Let~$X$ be a closed aspherical Kähler manifold of complex dimension $n\geq
  2$. Let~$\alpha$ be a holomorphic $1$-form on~$X$ with isolated zeroes and
  let $a=\Re \alpha$. Then there is a neighborhood~$U$ of the class of~$a$ in $H^1(X;
  \mathbf{R})$ such that for every~$b$ in $U\cap H^1(X; \mathbf{Q})$, the kernel of~$b$ is of
  type $\mathcal{F}_{n-1}$. If furthermore $X$~has nonzero
  Euler
  characteristic, then the kernel of~$b$ is not of type
  $\FP_n(\mathbf{Q})$.
\end{prop}

\begin{rema}
  Since $X$~is Kähler and closed, holomorphic $1$-forms are automatically harmonic
  and consequently closed. Furthermore, from the Hodge decomposition, the dimension of the
  space of holomorphic $1$-forms is half the first Betti number. Hence the
  assumption on~$X$ is of topological flavor.
\end{rema}

A deformation argument \parencite[section 6.2]{isenrich-py-exotic} shows that
the class of~$a$ can be
represented by a Morse $1$-form (i.e.\ locally the differential of a Morse
function) all of whose critical points have index equal to~$n$. This property
will hold in a neighborhood~$U$ of the class of~$a$ in $H^1(X; \mathbf{R})$
\parencite[proposition 8.1]{isenrich-py-exotic}. Let~$b$ be a rational form in
the open set~$U$ and choose~$\beta$ a differential form representing~$b$.

The universal cover~$\widetilde{X}$ of~$X$ is a contractible manifold and the
lift of~$\beta$ is  the differential of a function $\widetilde{X}\to
\mathbf{R}$.  This function descends to a function~$f\colon X_0\to
\mathbf{R}$, where 
$X_0=\widetilde{X}/ \ker b$ is the cover associated with~$b$.
The space~$X_0$ is aspherical with fundamental group equal to
$\ker b$, thus the finiteness properties of $\ker b$ can be determined
from~$X_0$ or from spaces homotopically equivalent to~$X_0$.
The function~$f$
is proper and has isolated singularities all of index~$n$. Therefore
Morse--Lefschetz theory implies that $X_0$~has the homotopy type of a compact
manifold (a regular fiber of~$f$) with infinitely many $n$-cells attached (as
soon as the form~$\alpha$ has at least one zero, which is ensured by the
assumption on the Euler characteristic). This model for the classifying space
of the group $\ker b$ implies that $\ker b$ is indeed of
type~$\mathcal{F}_{n-1}$. Using a long exact sequence due to
\textcite{milnor-cyclic} associated with the cyclic covering $
X_0 \to X$, \textcite[section~3.2]{isenrich-martelli-py} show that $\ker b$ is not
of type $\FP_n(\mathbf{Q})$.

\subsection{Maps to tori}
\label{sec:maps-tori}

Holomorphic forms on tori  are easily understood and never vanish (unless
zero). A way of obtaining holomorphic $1$-forms with isolated zeroes will be
by pulling back forms on tori.

\begin{prop}[{\cite{simpson}, cf.\ \cite[propositions 14 and 18]{llosa-pi-inventiones}}]
  Let~$X$ be a compact complex manifold and~$A$ be a complex torus. Let
  $\psi\colon X\to A$ be a holomorphic and finite-to-one map. There is then a meager
  set~$F$ in $H^0(A; \Omega^{1}_{A} )$ (i.e.\ $F$~is a countable union of
  closed nowhere dense subsets) such that, for every~$\beta$ in $H^0(A,
  \Omega^{1}_{A})\smallsetminus F$, the holomorphic $1$-form $\psi^* \beta$
  has isolated zeroes. 
\end{prop}

The argument goes as follows. Let~$Z$ be a connected component of the zeroes
of~$\alpha$ then $\beta$~must be zero on the subtorus generated by the image
of~$Z$ by~$\psi$. Since there are countably many nontrivial subtori, adjusting~$F$
appropriately, we can conclude that $\psi(Z)$ is a point and thus $Z$~is as
well a point since $\psi$~is finite-to-one.

\subsection{The Albanese map}
\label{sec:albanese-map}

Let~$X$ be a connected Kähler manifold.
Every complex differential $1$-form~$\alpha$ (and in particular every holomorphic
differential $1$-form) can be integrated along a path~$\gamma$ in~$X$ and the
resulting complex number $\int_\gamma \alpha$ depends only on the homotopy
class of~$\gamma$ (relative to the endpoints if any). There is thus a well
defined map from $H_1(X;\mathbf{Z})$ to  the dual space $H^0(X; \Omega^{0}_{X})^*$ whose image
is a lattice~$\Lambda$ in  $H^0(X; \Omega^{0}_{X})^*$. The
quotient of  $H^0(X; \Omega^{0}_{X})^*$ by~$\Lambda$ is called the
\emph{Albanese variety} of~$X$ and denoted by~$A(X)$.

Fixing a base point~$x_0$ in~$X$, we get a holomorphic map $a_X\colon X \to
A(X)$ called the \emph{Albanese map} as follows. For~$x$ in~$X$ choose a path
$\gamma_x$ starting from~$x_0$ and ending at~$x$ and set $a_X(x)$ to be the
class in $A(X)$ of the linear form $H^0(X, \Omega^{1}_{X})\to \mathbf{C}\mid
\alpha \mapsto \int_{\gamma_x} \alpha$; $a_X(x)$~does not depend on the choice
of~$\gamma_x$ precisely because of the quotient by the lattice~$\Lambda$.

The differential of the Albanese map is well understood (cf.\ lemma~23
in~\cite{llosa-pi-inventiones}) and this is one input for the following statement. 
\begin{theo}[{\cite[corollary 4.7]{eyssidieux2018orbifold}}]
  \label{theo:albanese-lattice-eyssidieux}
  Let~$\Gamma$ be an arithmetic lattice in $\mathrm{PU}(n,1)$ with positive first Betti
  number. There is then a finite index subgroup~$\Gamma_0$ of~$\Gamma$ such
  that the Albanese map for $X=B/\Gamma_0$ ($B$~being the unit ball
  in~$\mathbf{C}^n$) is an immersion and is thus finite-to-one.
\end{theo}
We refer to \cite{llosa-pi-inventiones} (theorem~24) for a proof.

\subsection{Lattices of the simplest type}
\label{sec:latt-simpl-type}

We now explain that there are indeed lattices in
$\mathrm{PU}(n,1)$ satisfying the hypothesis of
theorem~\ref{theo:albanese-lattice-eyssidieux}. This discussion is borrowed
from \cite{llosa-pi-inventiones} (section~3.2).

Let $F\subset \mathbf{R}$ be a totally real number field (for example
$F=\mathbf{Q}[\sqrt{2}]$) and let $E\subset \mathbf{C}$ be a purely imaginary
quadratic extension of~$F$ (for example $E=\mathbf{Q}[\sqrt{2},i]$). On
$V=E^{n+1}$ choose an Hermitian form of signature $(n,1)$ all of whose other
Galois transforms have definite signature (for example $(z_1, \dots,
z_{n+1})\mapsto z_1 \bar{z}_1 + \cdots + z_n \bar{z}_n -\sqrt{2}
z_{n+1}\bar{z}_{n+1}$). The group $\mathrm{U}(H, O_E)$ of $H$-Hermitian
matrices with coefficients in the ring~$O_E$ of integers of~$E$ (in the
example $O_E$ is $\mathbf{Z}[\sqrt{2}, i]$) is naturally a lattice in
$\mathrm{PU}(n,1)$. It is cocompact when $F$~is different from~$\mathbf{Q}$
and theorem~1 of \textcite{kazhdan} states that $\mathrm{U}(H, O_E)$ admits a finite
index subgroup with positive first Betti number. Since the quotient of the
unit ball in~$\mathbf{C}^n$ by a discrete cocompact subgroup has nonzero Euler
characteristic, the proposition~\ref{prop:forms-with-isolated} can be applied.

\begin{rema}
  \citeauthor{llosa-pi-inventiones} showed in fact the existence of infinitely
  many, pairwise not commensurable, word hyperbolic groups admitting subgroups
  of type~$\mathcal{F}_n$ and not of type~$\mathcal{F}_{n+1}$.
\end{rema}

\section{Construction from right-angled polytopes}
\label{sec:constr-from-right}

Hereafter the article \citeauthor{italiano-martelli-migliorini5}
 (\citeyear*{italiano-martelli-migliorini5}) will be
mentioned as~\cite{italiano-martelli-migliorini5} and the article \citeauthor{italiano2022hyperbolic}
(\citeyear*{italiano2022hyperbolic}) will be
mentioned as~\cite{italiano2022hyperbolic}.

The approach developed by~\cite{italiano-martelli-migliorini5} for the proof of
theorem~\ref{theo:les_italiens} is more combinatorial in nature. It starts by
constructing a finite volume hyperbolic $5$-manifold from a right-angled
polytope in the hyperbolic space~$\mathbb{H}_5$, then a fibration $f\colon
M\to S^1$ is constructed using a natural cubulation of the manifold~$M$. In
order to produce a compact object (and hence a word hyperbolic group) one needs to
cap the boundary components of~$M$ to obtained a metric space~$M^\vee$; this can be done maintaining the
negatively curved metric on~$M^\vee$ and an extension of the fibration exists.

\subsection{The polytope and the manifold}
\label{sec:polytope-manifold}

The chosen model for hyperbolic space~$\mathbb{H}_5$ is the Klein model:
the unit ball
in~$\mathbf{R}^5$ with geodesic given by Euclidean segments. 

The polytope~$P_5$ in~$\mathbb{H}_5$ is described as the intersections of the
half-spaces
\[ \underline{\varepsilon} \cdot \underline{x} = \sum_{i=1}^{5} \varepsilon_i
  x_i \leq 1, \quad \underline{x}\in \mathbb{H}_5\]
where $\underline{\varepsilon}$ varies in the subgroup of~$\{\pm 1\}^5$
defined by $\prod \varepsilon_i=1$, i.e.\ an even number of
the~$\varepsilon_i$ are equal to~$-1$. We refer to \cite{italiano-martelli-migliorini5} (section~1.1) for a
 complete description, and to \cite{italiano2022hyperbolic} for further details on that polytope as well
as a related series of right-angled polytopes in dimensions $3, \dots,
8$. These polytopes were previously studied by \textcite{potyagailo-vinberg}
who explained them starting from
certain hyperbolic simplices. They are related (by duality) to a series of
semiregular polytopes discovered by \textcite{gosset}.

The polytope~$P_5$ has finite volume, is right-angled, and has $16$~facets given by the
hyperplanes 
where equality is achieved in the equation above. It has a big group of
symmetries: the permutation of coordinates as well as the coordinate-wise
pluttifkation by $\underline{\varepsilon}$ (cf.\ \cite{brindacier}, for
this classical operation);
this produces a group of symmetries of type~$D_4$ and of order $2^4\times
5!=1920$. The hyperbolic reflections through the $16$~hyperplanes
bounding~$P_5$ generate a discrete subgroup~$\Gamma$ of $\mathrm{Isom}(\mathbb{H}_5)$
that is known to be isomorphic to the congruence two subgroup of the group of
integral matrices in the Lie group $\mathrm{O}(5,1)$ (see \textcite{ratcliffe-tschantz-5fld} who also
give a description of~$P_5$ in the hyperboloid model of hyperbolic space).

The group~$\Gamma$ is in fact a right-angled Coxeter group whose generating
system is given by the family $\{ r_F\}_F$ of reflections through the
facets~$F$ of~$P_5$; the relations, besides $r_{F}^{2}=e$, being $r_F r_G=r_G
r_F$ each time that two facets~$F$ and~$G$ intersect.

Each facet of~$P_5$ is adjacent to~$10$ other facets and this gives an adjacency
graph with $80$~edges (and $16$~vertices corresponding to the facets) that
controls the presentation of~$\Gamma$ and that can be nicely represented in
the plane (cf.\ figure~1 in~\cite{italiano-martelli-migliorini5}). 

By general properties of Coxeter groups, the torsion elements of~$\Gamma$ are
those conjugate to $r_{F_1} r_{F_2}\cdots r_{F_n}$ where $F_1, \dots, F_n$ are
facets of~$P_5$ that pairwise intersect. This happens only when $n\leq 5$ and
the $n$-tuples of facets satisfying this condition are completely determined by
the adjacency graph.

There is a natural homomorphism $\Gamma \to (\mathbf{Z}/2\mathbf{Z})^{16}$
whose kernel is torsion free and hence produces a hyperbolic, complete,
finite-volume manifold. To produce ``smaller'' manifolds,
\textcite{italiano-martelli-migliorini5, italiano2022hyperbolic}  use other
homomorphisms from~$\Gamma$ to $(\mathbf{Z}/2\mathbf{Z})^c$ (where $c$~is an
integer) that are in fact given combinatorially by a map from the facets
of~$P_5$ to $\{1, \dots, c\}$; the homomorphism $\Gamma\to
(\mathbf{Z}/2\mathbf{Z})^c$ is then uniquely determined by assigning to~$r_F$
the $i$-th basis element~$e_i$ of $(\mathbf{Z}/2\mathbf{Z})^c$ if $i$~is the integer
corresponding to~$F$. A necessary and sufficient condition for the kernel to
be torsion-free is then that adjacent facets are sent to different integers
under the mapping $F\mapsto i$. In the terminology of
\cite{italiano-martelli-migliorini5}, \cite{italiano2022hyperbolic} and others, the
mappings from the facets to $\{1, \dots, c\}$ are called colorings; they also
construct a coloring satisfying the above condition with $c=8$ (cf.\ figure~3
in \cite{italiano-martelli-migliorini5}).

The produced manifold~$M$ is hence made of $2^8$ copies of~$P_5$ which we
label $P_\lambda$ ($\lambda \in (\mathbf{Z}/2\mathbf{Z})^8$); along a facet~$F$
of~$P_\lambda$ is glued $P_{\lambda+e_i}$ where again $i$~is the ``color''
corresponding to the facet~$F$.

\subsection{The cusps of~$M$}
\label{sec:cusps-m}

The polytope~$P_5$ (that may be considered also as a hyperbolic orbifold)
has~$10$ cusps corresponding to the points in $\partial_\infty
\mathbb{H}_5\subset \mathbf{R}^5$ all of whose coordinates but one are equal
to~$0$, i.e.\ the points $(\pm 1,0,0,0,0),\dots, (0,0,0,0,\pm 1)$. They lift
to cusps in~$M$ and those lifts are analyzed in \cite{italiano-martelli-migliorini5} (section~1.4). The
coloring is not symmetric and the different lifts are not pairwise
isomorphic. The preimage of the cusp in~$P_5$ corresponding to one of the points
$(\pm 1,0,0,0,0),\dots, (0,0,0,\pm 1,0)$ is $1$~cusp in~$M$ and is named a
large cusp in \cite{italiano-martelli-migliorini5}. The preimage of a cusp in~$P_5$ corresponding to one of the
points $(0,0,0,0,\pm 1)$ consists of $2^4=16$ cusps in~$M$ and those cusps are called small in \cite{italiano-martelli-migliorini5}. There are thus $8+32=40$ cusps in~$M$.

The cusps in~$M$ naturally inherit a tessellation from the tessellation
of~$M$; the tiles of the cusps are $[0,1]^4\times \mathbf{R}_{\geq 0}$, the
product of the $4$-cube and the half line. The large cusps are divided in
$2^8=4^4$ tiles and are naturally isomorphic to
$(\mathbf{R}/4\mathbf{Z})^4\times \mathbf{R}_{\geq 0}$ with its ``natural''
tessellation coming from the tessellation of $\mathbf{R}^4\times
\mathbf{R}_{\geq 0}$ by translates of $[0,1]^4\times \mathbf{R}_{\geq 0}$. 
The small cusps are divided in $2^4$ tiles and are naturally isomorphic to
$(\mathbf{R}/2\mathbf{Z})^4\times \mathbf{R}_{\geq 0}$.

\subsection{The cubulation of~$M$}
\label{sec:cubulation-m}

The tessellation of~$M$ dual to the previous tessellation induces in fact a
tessellation of the compact manifold~$M^\vee$ obtained from~$M$ by removing
the cusps (therefore $M^\vee$~has $40$~toroidal boundary components). The
precise description starts in fact from the tessellation of~$M^\vee$ by copies
of~$P^\vee$, the polytope obtained from~$P_5$ by removing its cusps.  Taking the barycentric
subdivision of~$P^\vee$ produces then another tessellation of~$M^\vee$.  The
final tessellation is the one whose maximal polytopes are the stars, in this
intermediate subdivided tessellation, of
the vertices belonging to the
 tessellation of~$M^\vee$ by copies of~$P^\vee$ (we refer to
\cite{italiano-martelli-migliorini5}, section~1.5, for the precise construction). Since
$P_5$~is a right-angled polytope, this new tessellation is composed of cubes,
the vertices of this cubulation are in one-to-one correspondence with the
copies of~$P_5$ composing~$M$ (hence with $(\mathbf{Z}/2\mathbf{Z})^8$). The
edges of this cubulation are in one-to-one correspondence with the facets of
the original tessellation of~$M$ by copies of~$P_5$.

\subsection{The Bestvina--Brady Morse theory}
\label{sec:bestvina-brady-morse}

We give here a very quick sketch of this variant of Morse theory in the
piecewise linear setting. Let~$X$ be a topological space composed of copies of convex
polyhedra (in some finite dimensional real vector space) glued together via
affine maps; $X$~is called an affine cell complex. A function $f\colon X\to
\mathbf{R}$ will be said \emph{piecewise linear Morse} if (1) it is affine in restriction to
every cell (2) it is constant only on restriction to the $0$-dimensional cells,
and (3) the image of the $0$-skeleton is discrete \parencite[definition~2.2]{bestvina-brady}.

The link at a vertex~$x$ of~$C$ is defined to be the space made of the cells
containing~$x$ glued together along subcells containing~$x$.
When a piecewise linear Morse function~$f$ is given, the ascending link
$\mathrm{lk}_{\uparrow}(x,f)$ is the subspace of the link of~$X$ at~$x$ made
of the cells where $f$~attains its minimum at~$x$. Similarly the descending
link $\mathrm{lk}_{\downarrow}(x,f)$ is defined.

The topological changes of the sublevel sets $f^{-1}((-\infty, t])$ happen
only at vertices and are controlled by the ascending and descending links. The statement
that will be used below is that $X$~is homotopically equivalent to a fiber
of~$f$ when all the ascending and descending links are contractible.

\subsection{Affine maps from~$M$ to the circle}
\label{sec:affine-maps-from}

In order to apply the above Morse theory, we need to construct a piecewise linear map
from~$M$, or rather from $M^\vee$ with its affine cell complex structure
inherited from its cubulation, to the circle $S^1= \mathbf{R}/\mathbf{Z}$ with
its natural affine structure. The lift of this map is a map from the universal cover
of~$M^\vee$ to~$\mathbf{R}$ and the ascending and descending links of the lift
are exactly those of the initial map $M^\vee \to S^1$.

The map~$f$ from~$M^\vee$ to~$S^1$ will send all
the vertices of the cubulation to~$0$. There are then two possibilities for
its restriction to a given edge (a little abusively, once an identification of
the edge with the interval $[0,1]$ is given, the two possibilities are
the maps $x\mapsto x$ and $x \mapsto -x$). Once such choices on the edges are
given, this produces an piecewise linear map from the $1$-skeleton of~$M^\vee$
to~$S^1$. It is possible to extend it to~$M^\vee$ if and only if, for every
square $C\subset M^\vee$, the map $\partial C\to S^1$, deduced from the
inclusion of~$\partial C$ in the $1$-skeleton of~$M^\vee$, admits an affine
extension to~$C$, i.e.\ if this map is given by $(x,y)\mapsto x+y$ after
appropriately identifying~$C$ with $[0,1]^2$.

In the present situation, as edges of the cubulation of~$M^\vee$ are in
one-to-one correspondence with the facets of the tessellation of~$M$, we need
to co-orient the facets in~$M$. In fact, all the facets of all the copies~$P_\lambda$ of
the polytope $P_5$ will be co-oriented. For a facet~$F$ of~$P_\lambda$ we have
thus two possible co-orientations, either inward or outward. The
facet~$F$ also belong to~$P_{\lambda'}$ with $\lambda'-\lambda=e_i$ ($i$~being
the color of~$F$) and, at the very least, the co-orientations of~$F$ in
$P_\lambda$ and $P_{\lambda'}$ must be opposite.

The co-orientations of the facets are here determined algorithmically by
deciding the co-orientations of the facets of~$P_{\lambda+e_i}$ from the
knowledge of the co-orientations of the facets of~$P_\lambda$. At least,
all the facets whose color is~$i$ must have their co-orientations reverted. It
was observed in \cite{italiano2022hyperbolic} (proposition 12.(1)) that it is
necessary (and sufficient) to have two facets of~$P_\lambda$ with the same
color and opposite co-orientations in order to produce a nonzero homomorphism
from $\pi_1(M)$ to~$\mathbf{Z}$. 
Furthermore, reverting the co-orientation only for facets of the same color
cannot lead to maps satisfying all the desired properties (cf.\
\cite{italiano-martelli-migliorini5}, end of section~1.8; the last condition
needed on the co-orientations will be mentioned in
section~\ref{sec:restriction-cusps} below). 
One must then revert more co-orientations and this is done via
the following procedure: an equivalence relation~$\mathcal{R}$ on the colors $\{1, \dots, c\}$ (here
$c=8$) is chosen and, when going from~$P_{\lambda}$ to~$P_{\lambda+e_i}$ (for
all~$\lambda$ and for all~$i$), one reverts the co-orientation of all the
facets whose color~$j$ is equivalent to~$i$. This is expressed in terms of
partitions of $\{1,\dots, c\}$ in \cite{italiano-martelli-migliorini5}, \cite{italiano2022hyperbolic} and other references. This
procedure was developed by \textcite{jankiewicz-norin-wise} who introduced a specific vocabulary for it
(status, state, moves, game) that we did not reproduce
here. 

The chosen equivalence relation in~\cite{italiano-martelli-migliorini5} is the
following: $i\mathcal{R}j$ if and only if
$i=j\mod 4$. It is not difficult to determine directly the co-orientation
on the facets of~$P_\lambda$ in terms of a fixed co-orientation on the facets
of~$P_5$ and an homomorphism $(\mathbf{Z}/2\mathbf{Z})^8
\to(\mathbf{Z}/2\mathbf{Z})^4$ \parencite[section~1.6]{italiano-martelli-migliorini5}.
 This gives rise to a piecewise linear Morse map from~$M^\vee$ to~$S^1$.

\subsection{The restriction to the cusps}
\label{sec:restriction-cusps}

 We discuss here the restrictions of the piecewise linear map to the boundary components
of~$M^\vee$. This gives maps from the cubulated tori
$(\mathbf{R}/2\mathbf{Z})^4$ or $(\mathbf{R}/4\mathbf{Z})^4$ to the circle. In
order to produce later relevant compact objects, we will need to cap off these
boundary tori and to extend non-trivially the map to the circle. This will be
possible when the restrictions of the piecewise linear map to the tori are homotopically nonzero. 

For the coloring given in the previous subsection, the homotopy classes of the
restrictions to boundary tori are calculated in~\cite{italiano-martelli-migliorini5} (proposition~14); 
the restriction to a large cusp is homotopic to the projection $(S^1)^4\to
S^1$ on a factor; the restriction to a small cusp
 is homotopic to the summation map $(\mathbf{R}/\mathbf{Z})^4\to
\mathbf{R}/\mathbf{Z}\mid (x_1,x_2,x_3,x_4)\mapsto x_1+x_2+x_3+x_4$.

\subsection{A complication}
\label{sec:complication}

As explained in section~\ref{sec:affine-maps-from}, 
the above equivalence relation give rise to a piecewise linear map from the
$1$-skeleton of~$M^\vee$ to the circle. However this map cannot extend to a
piecewise linear map on~$M^\vee$ since some
squares in~$M^\vee$ do not satisfy the condition stated in section~\ref{sec:affine-maps-from}. Even worse: no equivalence relation on
$\{1, \dots, 8\}$ satisfies all the wanted properties, i.e.\ the conditions on
squares and the conditions on boundary components
(see~\cite{italiano-martelli-migliorini5}, section~1.8).

This issue can be circumvented as follows. \cite{italiano-martelli-migliorini5} analyzes the ``bad'' squares~$C$ and
proves that they always appear in $5$-cubes $C\times D$ ($D$~is hence a
$3$-cube) such that every parallel copy $C\times \{x\}$ ($x$~vertex of~$D$)
is bad (\cite{italiano-martelli-migliorini5}, proposition~12). This enables to further subdivide these
$5$-cubes, and their subcubes,
into prisms $T\times D$ ($T$ being one of the $4$~triangles obtained by
cutting~$C$ along its diagonals) and leaving the other $5$-cubes unaffected. On
this new tessellation of $M^\vee$ everything goes well: the map extends, the
links are contractible (\cite{italiano-martelli-migliorini5}, theorem~13) and the restrictions to boundary
components have the form mentioned previously (in fact proposition~14 of
\cite{italiano-martelli-migliorini5} mentioned above
calculates with this new tessellation).

\subsection{Capping off the boundary components}
\label{sec:capp-bound-comp}

This process is nicely explained in different places, for example in \textcite[section~2.1]{isenrich-martelli-py}. For each
boundary component~$T$ of $M^\vee$, one glues to~$M^\vee$ along~$T$ the space
$(T\times [0,1])/{\sim}$ where $\sim$~is the equivalence relation whose
classes are $\{t,s\}$ ($t\in T$, $s<1$) and $S\times \{1\}\subset T\times [0,1]$ ($S\subset T$ is a fiber
of the map $T\to S^1$ obtained by restricting the map $M^\vee\to S^1$). The map $M^\vee\to S^1$
obviously extends to the reunion.

The hyperbolic metric on~$M$ induces a flat metric on~$T$. \textcite[theorem~2.7]{fujiwara-manning} 
ensure that, when $T$~is ``big enough'' (the shortest loop has length at least~$2\pi$), then the obtained space carries a metric that is locally
$\mathrm{CAT}(-1)$. The $2\pi$-condition for all the boundary tori can be achieved up to taking a
finite cover~$N$ of~$M$. Capping off all the boundary tori gives a compact
pseudo-manifold $N^\dagger$ with a locally $\mathrm{CAT}(-1)$ metric and a map
$N^\dagger\to S^1$ which is a fibration. The fiber $F^\dagger$ of this map is
also a pseudo-manifold which can be assumed to be connected
(cf.~\cite{italiano-martelli-migliorini5}, remark~16).

\subsection{Asphericity and non-hyperbolicity}
\label{sec:asph-non-hyperb}

The fundamental group $\pi_1( N^\dagger)$ is then word hyperbolic and contains the
fundamental group $\pi_1( F^\dagger)$ as a normal subgroup. Application of the
Bestvina--Brady Morse theory implies that the universal cover of~$F^\dagger$
is homotopically equivalent to the universal cover of $N^\dagger$ and is hence
contractible: the group $\pi_1(F^\dagger)$ has type~$\mathcal{F}$.

The fact that $\pi_1(F^\dagger)$ is not word hyperbolic is shown as follows (\cite{italiano-martelli-migliorini5}
section~3). Its outer automorphism group is infinite
(\cite{italiano-martelli-migliorini5}, proposition~23), if it were
word hyperbolic, it would split over a cyclic group
\parencite[corollary~1.3]{bestvina-feighn}, but a
Mayer--Vietoris argument shows that this is not possible
(\cite{italiano-martelli-migliorini5}, proposition~24).

\begin{rema}
  Using the same arguments, every fibration $M\to S^1$ of a closed
  hyperbolic manifold of odd dimension $\geq 5$ will give nonhyperbolic
  subgroups of type~$\mathcal{F}$ (the fundamental group of a fiber) in the
  fundamental group of~$M$.
\end{rema}


\printbibliography

\end{document}
